\documentclass[11pt, a4paper,twoside]{article}

\usepackage{amssymb}
\usepackage{amsmath}
\usepackage{fancyhdr}
\usepackage{amsthm}
\usepackage{afterpage}  

\newtheorem{theorem}{Theorem}

\newtheorem{corollary}{Corollary}
\newtheorem{conjecture}{Conjecture}

\def\tref#1{Theorem~$\ref{#1}$}
\def\Tref#1{Table~\ref{#1}}

\def\cref#1{Conjecture~$\ref{#1}$}
\def\Cref#1{Corollary~\ref{#1}}
\def\eref#1{(\ref{#1})}
\def\sref#1{\S\ref{#1}}

\def\Z{\mathbb{Z}}

\def\mod{{\ \rm mod\ }}
\def\id{\varepsilon}

\def\dfrac#1#2{\lower0.15ex\hbox{\large$\frac{#1}{#2}$}}

\begin{document}

\afterpage{\rhead[]{\thepage} \chead[\small I. M. Wanless]{\small  Transversals in Latin Squares} \lhead[\thepage]{} }

\begin{center}

\vspace*{2pt}

{\Large \textbf{Transversals in Latin Squares}}\\[36pt]

{\large \textsf{\emph{Ian M. Wanless}}}\\[36pt]

\textbf{Abstract}

\end{center}

{\footnotesize 
A latin square of order $n$ is an {$n\times n$} array of $n$ symbols
in which each symbol occurs exactly once in each row and column. A
transversal of such a square is a set of $n$ entries such that no two
entries share the same row, column or symbol. Transversals are closely
related to the notions of complete mappings and orthomorphisms in
(quasi)groups, and are fundamental to the concept of mutually orthogonal
latin squares.

\footnotesize 
Here we provide a brief survey of the literature on transversals. We
cover (1)~existence and enumeration results, (2)~generalisations of
transversals including partial transversals and plexes, (3)~the
special case when the latin square is a group table, (4)~a connection
with covering radii of sets of permutations.
The survey includes a number of conjectures and open problems.

}



\section*{\centerline{1. Introduction}}%\label{s:intro}

A \emph{latin square} of order $n$ is an {$n\times n$} array of $n$
symbols in which each symbol occurs exactly once in each row and
in each column. By a \emph{diagonal} of such a square we mean a set
of entries which contains exactly one representative of each row and 
column. A \emph{transversal} is a diagonal in which no symbol is repeated.

Historically, interest in transversals arose from the study of
orthogonal latin squares. A pair of latin squares $A=[a_{ij}]$ and
$B=[b_{ij}]$ of order $n$ are said to be \emph{orthogonal mates} if the
$n^2$ ordered pairs $(a_{ij},b_{ij})$ are distinct. It is simple to
see that if we look at all $n$ occurrences of a given symbol in $B$,
the corresponding positions in $A$ must form a transversal. Indeed,

\begin{theorem}\label{t:transMOLS}
A latin square has an orthogonal mate iff it has a decomposition into
disjoint transversals.
\end{theorem}

For example, below there are two orthogonal latin squares of order 8.
Subscripted letters are used to mark the transversals of the left hand
square which correspond to the positions of each symbol in its orthogonal
mate (the right hand square).
\begin{equation}\label{e:ols8}
\begin{matrix}
1_a&2_b&3_c&4_d&5_e&6_f&7_g&8_h\\ 7_b&8_a&5_d&6_c&2_f&4_e&1_h&3_g\\
2_c&1_d&6_a&3_b&4_g&5_h&8_e&7_f\\ 8_d&7_c&4_b&5_a&6_h&2_g&3_f&1_e\\
4_f&3_e&1_g&2_h&7_a&8_b&5_c&6_d\\ 6_e&5_f&7_h&8_g&1_b&3_a&2_d&4_c\\
3_h&6_g&2_e&1_f&8_c&7_d&4_a&5_b\\ 5_g&4_h&8_f&7_e&3_d&1_c&6_b&2_a
\end{matrix}
\qquad\begin{matrix}
a&b&c&d&e&f&g&h\\ b&a&d&c&f&e&h&g\\
c&d&a&b&g&h&e&f\\ d&c&b&a&h&g&f&e\\
f&e&g&h&a&b&c&d\\ e&f&h&g&b&a&d&c\\
h&g&e&f&c&d&a&b\\ g&h&f&e&d&c&b&a
\end{matrix}
\end{equation}
%% Second square has 128 transversals (first only has 8)

More generally, there is interest in sets of {\em mutually orthogonal
latin squares} (MOLS), that is, sets of latin squares in which each
pair is orthogonal in the above sense.  The literature on MOLS is vast
(start with \cite{DKI,DKII,LM}) and provides ample justification for
an interest in transversals. Subsequent investigations have ranged far
beyond the initial justification of \tref{t:transMOLS} and have proved
that transversals are interesting objects in their own right.  Despite
this, a number of basic questions about their properties remain
unresolved, as will become obvious in the subsequent pages.

Orthogonal latin squares exist for all orders $n\not\in\{2,6\}$.  For
$n=6$ there is no pair of orthogonal squares, but we can get
close. Finney~\cite{BFinnI} gives the following example 
which contains 4 disjoint transversals indicated by the subscripts
$a,b,c$ and $d$.
$$\begin{matrix}
1_a\hfill&2\hfill&3_b\hfill&4_c\hfill&5\hfill&6_d\hfill\cr
2_c\hfill&1_d\hfill&6\hfill&5_b\hfill&4_a\hfill&3\hfill\cr
3\hfill&4_b\hfill&1\hfill&2_d\hfill&6_c\hfill&5_a\hfill\cr
4\hfill&6_a\hfill&5_c\hfill&1\hfill&3_d\hfill&2_b\hfill\cr
5_d\hfill&3_c\hfill&2_a\hfill&6\hfill&1_b\hfill&4\hfill\cr
6_b\hfill&5\hfill&4_d\hfill&3_a\hfill&2\hfill&1_c\hfill\cr 
\end{matrix}$$ 

%That is, each main class representative $M$ is scored by the maximum
%$m$ such that $M$ has a $(1^m,n-m)$ partition.  Recall that if $m=n$
%then $M$ has a 1-partition and hence an orthogonal mate.

\begin{table}[htb]
\centering
\begin{tabular}{|c|ccccc|}       
\hline
$m$&$n=4$&5&6&7&8\\
\hline
0&1&0&6&0&33\\
1&0&1&0&1&0\\
2&0&0&2&5&7\\
3&-&0&0&24&46\\
4&1&-&4&68&712\\
5&-&1&-&43&71330\\
6&-&-&0&-&209505\\
7&-&-&-&6&-\\
8&-&-&-&-&2024\\
\hline
Total&2&2&12&147&283657\\
\hline
\end{tabular}
\caption{\label{T:distrans}Number $m$ of disjoint transversals 
in latin squares of order $n\leq8$.}
\end{table}

\Tref{T:distrans} shows the squares of order $n$, for $4\leq n\leq 8$,
counted according to their maximum number $m$ of disjoint
transversals. The entries in the table are counts of main classes (A
{\em main class}, or {\em species} is an equivalence class of latin
squares each of which has essentially the same structure. See
\cite{DKI,LM} for the definition.)

Evidence such as that in \Tref{T:distrans} led van Rees \cite{vnRs} to
conjecture that, as $n\rightarrow\infty$, a vanishingly small
proportion of latin squares have orthogonal mates. However, the trend
seems to be quite the reverse (see \cite{WW}), although no rigorous way of
establishing this has yet been found.

A point that \Tref{T:distrans} raises is that some latin squares
have no transversals at all. We now look at some results in this regard.

A latin square of order $mq$ is said to be of {\it$q$-step type\/} if
it can be represented by a matrix of $q\times q$ blocks $A_{ij}$ as
follows
$$\begin{matrix}\label{eqstep}
A_{11}&A_{12}&\cdots&A_{1m}\cr
A_{21}&A_{22}&\cdots&A_{2m}\cr
\vdots&\vdots&\ddots&\vdots\cr
A_{m1}&A_{m2}&\cdots&A_{mm}\cr
\end{matrix}$$
where each block $A_{ij}$ is a latin subsquare of order $q$
and two blocks $A_{ij}$ and $A_{i'j'}$ contain the same symbols
iff $i+j\equiv i'+j'\mod m$. The following classical theorem is
due to Maillet \cite{Maillet}.

\begin{theorem}\label{Tqstep} 
Suppose that $q$ is odd and $m$ is even.
No $q$-step type latin square of order $mq$ possesses a transversal.
\end{theorem}

As we will see in \sref{s:groups}, this rules out many
group tables having transversals. In particular, 
no cyclic group of even order has a transversal. By contrast,
there is no known example of a latin square of odd order without
transversals. 

\begin{conjecture}\label{c:ryser}
Each latin square of odd order has at least one transversal.
\end{conjecture}

This conjecture is known to be true for $n\le9$ (see \sref{s:enum}).  It is
attributed to Ryser~\cite{hr67} and has been open for forty years. In
fact, Ryser's original conjecture was somewhat stronger: for every
latin square of order~$n$, the number of transversals is congruent to
$n\mod2$. In~\cite{kb90}, Balasubramanian proved the even case.

\begin{theorem}\label{t:balasub}
In any latin square of even order the number of transversals is even.
\end{theorem}

Despite this, it has been noted in \cite{CW} (and other places) that
there are many counterexamples of odd order to Ryser's original
conjecture. Hence the conjecture has now been weakened to
\cref{c:ryser} as stated.  One obstacle to proving this conjecture was
recently revealed in \cite{WW}.

\begin{theorem}\label{t:bachelor}
For every $n>3$ there exists a latin square of order $n$ which contains
an entry that is not included in any transversal.
\end{theorem}

Given \tref{t:transMOLS}, this latest theorem showed existence for all
$n>3$ of a latin square without an orthogonal mate. The same result was
obtained in \cite{EvansBach} without showing \tref{t:bachelor}.

\section*{\centerline{2. Partial transversals}}%\label{s:partial}

We have seen in \sref{s:intro} that not all latin squares have
transversals, which prompts the question of how close we can get to
finding a transversal in such cases. We define a {\em partial
transversal of length $k$} to be a set of $k$ entries, each selected
from different rows and columns of a latin square such that no two
entries contain the same symbol. Note that in some papers
(e.g.~\cite{ps82}) a partial transversal of length $k$ is defined
slightly differently to be a diagonal on which $k$ different symbols
appear.

Since not all squares of order $n$ have a partial transversal of
length $n$ (i.e.~a transversal), the best we can hope for is to find
one of length $n-1$.  The following conjecture has been attributed by
Brualdi (see \cite[p.103]{DKI}). 

\begin{conjecture}\label{c:brualdi}
Every latin square of order $n$ possesses a
partial transversal of length $n-1$. 
\end{conjecture}

A claimed proof of this conjecture by Derienko~\cite{id88}
contains a fatal error~\cite{CW}. Recently, a paper \cite{HuLi} has appeared in
the maths arXiv claiming a proof of \cref{c:brualdi}. 
Given the history of the problem such a claim should be treated cautiously,
at least until the paper has been refereed.\looseness-1

The best reliable result to date states that there must be a partial
transversal of length at least $n-O(\log^2n)$. This was shown by
Shor~\cite{ps82}, and the implicit constant in the `big \emph{O}' was
very marginally improved by Fu et al.~\cite{hf02}. Subsequently Hatami
and Shor \cite{ShorErratum} discovered an error in \cite{ps82}
(duplicated in \cite{hf02}) and corrected the constant to a higher
one. Nonetheless, the important thing remains that the bound is
$n-O(\log^2n)$. This improved on a number of earlier bounds including
$\frac23n+O(1)$ (Koksma \cite{Koksma}), $\frac34n+O(1)$ (Drake
\cite{Drake}) and $n-\surd n$ (Brouwer et al.~\cite{BVW} and
Woolbright \cite{Woolbright}).\looseness-1

%% There is one more in DKII -- Wang's Phd thesis showed 9n/11

Erd\H os and Spencer \cite{ES} showed that any $n\times n$ array in
which no entry occurs more than $(n-1)/16$ times has a transversal (in
the sense of a diagonal with $n$ different symbols on it).  It has
also been shown by Cameron and Wanless~\cite{CW} that every latin
square possesses a diagonal in which no symbol appears more than
twice.

\cref{c:brualdi} has been well known and open for decades. A much
simpler problem is to consider the shortest possible length of a
maximal partial transversal (maximal in the sense that it is contained
in no partial transversal of greater length). It is easy to see that
no partial transversal of length strictly less than $\frac12n$ can be
maximal, since there are not enough `used' symbols to fill the
submatrix formed by the `unused' rows and columns. However, for all
$n>4$, maximal partial transversals of length
$\big\lceil\frac12n\big\rceil$ can easily be constructed using a
square of order $n$ which contains a subsquare $S$ of order
$\big\lfloor\frac12n\big\rfloor$ and a partial transversal containing
the symbols of $S$ but not using any of the same rows or columns as $S$.

\section*{\centerline{3. Number of transversals}}%\label{s:enum}

In this section we consider the question of how many transversals a
latin square can have. We define $t(n)$ and $T(n)$ to be respectively
the minimum and maximum number of transversals among the latin squares
of order~$n$.

We have seen in \sref{s:intro} that some latin squares have no
transversals but it is not settled for which orders such latin squares
exist.  Thus for lower bounds on $t(n)$ we cannot do any better than
to observe that $t(n)\ge0$, with equality occurring at least when $n$
is even.  A related question, for which no work seems to have been
published, is to find an upper bound on $t(n)$ when $n$ is odd.

Turning to the maximum number of transversals, it should be clear that
$T(n)\le n!$ since there are only $n!$ different diagonals. An
exponential improvement on this trivial bound was obtained by McKay et
al.~\cite{MMW}:

\begin{theorem}\label{t:Tupper}
For $n \geq 5$,
\begin{equation*}
	15^{n/5} \leq T(n) \leq c^n\sqrt{n}\:n!  
\end{equation*}
where $c = \sqrt{\frac{3 - \sqrt{3}}{6}} \, e^{\sqrt{3}/6} \approx 0.61354$.
\end{theorem}

The lower bound in \tref{t:Tupper} is very simple and would not be too
difficult to improve. The upper bound took considerably more work,
although it too is probably far from the truth. 

In the same paper the authors reported the results of an exhaustive computation
of the transversals in latin squares of orders up to and including 9.
Table~\ref{T:small} lists the minimum and maximum number of
transversals over all latin squares of order $n$ for {$n \leq 9$}, and
the mean and standard deviation to 2 decimal places.

\begin{table}[htb]
\centering
\begin{tabular}{|l|cccc|}       \hline
$n$	& $t(n)$& Mean  & Std Dev	& $T(n)$\\ \hline
2 	& 0   	& 0     & 0     	& 0		\\      
3 	& 3   	& 3     & 0     	& 3		\\
4 	& 0	& 2     & 3.46 		& 8		\\
5 	& 3   	& 4.29  & 3.71          & 15		\\
6 	& 0   	& 6.86  & 5.19   	& 32		\\
7 	& 3   	& 20.41 & 6.00		& 133		\\  
8 	& 0   	& 61.05 & 8.66		& 384	        \\
9	& 68  	& 214.11& 15.79  	& 2241  	\\ \hline
\end{tabular}
\caption{\label{T:small}Transversals in latin squares of order $n\leq9$.}
\end{table}
 
%There are only three main classes of order~9 which possess more than
%1000 transversals. They are $\Z_3\oplus\Z_3$, which has 2241, $\Z_9$
%which has 2025 transversals and a third square, $L$, which has 1620
%transversals. Like the two group based squares, $L$ is composed of
%9~disjoint order~3 subsquares and has no intercalates ($2\times2$
%Latin subsquares).

\Tref{T:small} confirms \cref{c:ryser} for $n\leq9$.  The following
semisymmetric squares (see \cite{DKI} for a definition of
semisymmetric) are representatives of the unique main class with
$t(n)$ transversals for $n\in\{5,7,9\}$. In each case the largest
subsquares are shown in \textbf{bold}.
\begin{equation*}
	\begin{matrix}
	   \bf1&\bf2&\bf3&\bf4&\bf5\cr
	   \bf2&\bf1&4&5&3\cr
	   \bf3&5&\bf1&2&4\cr
	   \bf4&3&5&\bf1&2\cr
	   \bf5&4&2&3&\bf1\cr
	\end{matrix}\qquad
	\begin{matrix}
	  \bf3&\bf2&\bf1&5&4&7&6\cr
	  \bf2&\bf1&\bf3&6&7&4&5\cr
	  \bf1&\bf3&\bf2&7&6&5&4\cr
	  5&6&7&4&1&2&3\cr
	  4&7&6&1&5&3&2\cr
	  7&4&5&2&3&6&1\cr
	  6&5&4&3&2&1&7\cr
	\end{matrix}\qquad
	\begin{matrix}
		\bf2&\bf1&\bf3&6&7&8&9&5&4\cr
		\bf1&\bf3&\bf2&5&4&9&\bf6&\bf7&\bf8\cr
		\bf3&\bf2&\bf1&4&9&5&\bf7&\bf8&\bf6\cr
		9&5&4&3&2&1&\bf8&\bf6&\bf7\cr
		8&4&6&2&5&7&1&9&3\cr
		4&\bf7&\bf9&\bf8&3&6&5&1&2\cr
		5&\bf8&\bf7&\bf9&6&\bf2&\bf3&\bf4&1\cr
		6&\bf9&\bf8&\bf7&1&\bf4&\bf2&\bf3&5\cr
		7&6&5&1&8&\bf3&\bf4&\bf2&9\cr
	\end{matrix}
\end{equation*}

\begin{table}[htb]
	\centering
	\begin{tabular}{|c|cc|}       \hline
	$n$	& Lower Bound		& Upper Bound	\\ \hline
	10 	& 5504			& 75000		\\      
	11 	& 37851			& 528647	\\
	12 	& 198144		& 3965268	\\
	13 	& 1030367		& 32837805	\\
	14 	& 3477504		& 300019037	\\
	15 	& 36362925		& 2762962210	\\  
	16 	& 244744192		& 28218998328	\\
	17	& 1606008513		& 300502249052	\\ 
	18 	& 6434611200 		& 3410036886841	\\
	19 	& 87656896891		& 41327486367018\\ 
	20	& 697292390400	        & 512073756609248 \\
	21	& 5778121715415         & 6803898881738477 \\ \hline
	\end{tabular}
	\caption{\label{Tab:med}Bounds on $T(n)$ for $10\le n\le21$.}
\end{table}

In Table~\ref{Tab:med} we reproduce from \cite{MMW} bounds on $T(n)$
for {$10\leq n\le21$}. The upper bound is somewhat sharper than that
given by \tref{t:Tupper}, though proved by the same methods.  The
lower bound in each case is constructive and likely to be very close
to the true value.  When $n\not\equiv2\mod4$ the lower bound comes
from the group with the highest number of transversals (see
Table~\ref{Tab:groups}).  When $n\equiv2\mod4$ the lower bound comes
from a so-called turn-square, many of which were analysed in
\cite{MMW}.  A {\em turn-square} is obtained by starting with the
Cayley table of a group (typically a group of the form $\Z_2\oplus
\Z_m$ for some $m$) and ``turning'' some of the intercalates (that is,
replacing a subsquare of order 2 by the other possible subsquare on the same
symbols). For example,
\begin{equation}\label{e:parker}
    \begin{array}{ccccc|ccccc}
    	\bf5&\bf6&2&3&4   &	\bf0&\bf1&7&8&9 \\
	    \bf6&2&3&4&0  &    \bf1&7&8&9&5 \\
	    2&3&4&0&1 	  &  7&8&9&5&6 \\
	    3&4&0&1&2 	  &  8&9&5&6&7 \\
	    4&0&1&2&3 	  &  9&5&6&7&8 \\
\hline
    	\bf0&\bf1&7&8&9   & 	\bf5&\bf6&2&3&4 \\
	    \bf1&7&8&9&5  &   \bf6&2&3&4&0 \\
	    7&8&9&5&6 	  &  2&3&4&0&1 \\
	    8&9&5&6&7 	  &  3&4&0&1&2 \\
	    9&5&6&7&8 	  &  4&0&1&2&3 \\
    \end{array}
\end{equation}
achieves $5504$ transversals.  The `turned' entries have been marked
in \textbf{bold}.  The study of turn-squares was pioneered by Parker
(see \cite{bp93} and the references therein) in his unsuccessful quest
for a triple of MOLS of order 10.  He noticed that turn-squares often
have many more transversals than is typical for squares of their
order, and used this as a heuristic in the search for MOLS.

It is has long been suspected that $T(10)$ is achieved by
\eref{e:parker}.  This suspicion was strengthened by McKay et
al.~\cite{mmm05} who examined several billion squares of order 10,
including every square with a non-trivial symmetry, and found none had
more than $5504$ transversals.  Parker was indeed right that the
square \eref{e:parker} is rich in orthogonal mates (it has 12265168 of
them~\cite{mw04}, which is an order of magnitude greater than he
estimated). However, using the number of transversals as a heuristic
in searching for MOLS is not fail-safe. For example, the turn-square
of order 14 with the most transversals (namely, 3477504) does not have
any orthogonal mates \cite{MMW}. Meanwhile there are squares of order
$n$ with orthogonal mates but which possess only the bare minimum of
$n$ transversals (the left hand square in \eref{e:ols8} is one such).

Nevertheless, the number of transversals does provide a useful
invariant for squares of small orders where this number can be
computed in reasonable time (see, for example, \cite{krstdk96} and
\cite{iw02}).  It is straightforward to write a backtracking algorithm
to count transversals in latin squares of small order, though this
method currently becomes impractical if the order is much over 20. See
\cite{hhs}, \cite{hsc} for some algorithms and complexity theory
results on the problem of counting transversals.

It seems very difficult
to find theoretical estimates for the number of transversals (unless,
of course, that number is zero). This difficulty is so acute that
there are not even good estimates for $z_n$, the number of
transversals of the cyclic group of order $n$. 
Vardi~\cite{iv91} makes the following prediction:\looseness-1

\begin{conjecture}\label{c:vardi}
There exist real constants $0<c_1<c_2<1$ such that
\begin{equation*}
        c_1^nn! \leq z_n \leq c_2^nn!
\end{equation*}
for all odd $n\geq3$.
\end{conjecture}

Vardi makes this conjecture while considering a variation on the
toroidal {$n$-queens} problem. The toroidal $n$-queens problem is that
of determining in how many different ways $n$ non-attacking queens can
be placed on a toroidal {$n \times n$} chessboard. Vardi considered
the same problem using semiqueens in place of queens, where a
semiqueen is a piece which moves like a toroidal queen except that it
cannot travel on right-to-left diagonals. The solution to Vardi's
problem provides an upper bound on the toroidal $n$-queens
problem. The problem can be translated into one concerning latin
squares by noting that every configuration of $n$ non-attacking
semiqueens on a toroidal {$n \times n$} chessboard corresponds to a
transversal in a cyclic latin square $L$ of order $n$, where {$L_{ij}
\equiv i-j \mod n$}. Note that the toroidal $n$-queens problem is
equivalent to counting diagonals which simultaneously yield
transversals in $L$ and $L'$, where {$L_{ij}'=i+j \mod n$}.

As a corollary of \tref{t:Tupper} we can infer that the upper bound in
\cref{c:vardi} is true (asymptotically) with $c_2=0.614$. This also yields
an upper bound for the number of solutions to the toroidal $n$-queens
problem. \tref{t:Tupper} is valid for all latin squares, but
\cref{c:vardi} has also been attacked by methods which are specific to
the cyclic square.  Cooper and Kovalenko \cite{CK} first showed that
Vardi's upper bound is asymptotically true with $c_2=0.9153$, and this
was then improved to $c=1/\surd2\approx0.7071$ in \cite{Koval}.
Finding a lower bound of the form given in Conjecture~\ref{c:vardi} is
still an open problem. However, \cite{Cooper} and \cite{iv94} do give
some lower bounds, each of which applies only for some~$n$.  Cooper
et~al.~\cite{CGKN} estimated that perhaps the correct rate of growth
for $z_n$ is around $0.39^n\,n!$.

\section*{\centerline{4. Finite Groups}}%\label{s:groups}

By using the symbols of a latin square to index its rows and columns,
each latin square can be interpreted as the Cayley table of a
quasigroup.  In this section we consider the important special case
when that quasigroup is associative; in other words, it is a
group.

Much of the study of transversals in groups has been phrased in terms
of the equivalent concepts of complete mapping and
orthomorphisms. Mann \cite{Mann} introduced complete mappings for
groups, but their definition works just as well for quasigroups. It is
this: a permutation $\theta$ of the elements of a quasigroup
$(Q,\oplus)$ is a {\em complete mapping} if $\eta:Q\mapsto Q$ defined
by $\eta(x)=x\oplus\theta(x)$ is also a permutation.  The permutation
$\eta$ is known as an {\em orthomorphism} of $(Q,\oplus)$, following
terminology introduced in \cite{JDM}.  All of the results of this
paper could be rephrased in terms of complete mappings and/or
orthomorphisms because of our next observation.

\begin{theorem}
Let $(Q,\oplus)$ be a quasigroup and $L_Q$ its Cayley table.  Then
$\theta:Q\mapsto Q$ is a complete mapping iff we can locate a
transversal of $L_Q$ by selecting, in each row $x$, the entry in
column $\theta(x)$. Similarly, $\eta:Q\mapsto Q$ is an orthomorphism
iff we can locate a transversal of $L_Q$ by selecting, in each row
$x$, the entry containing symbol $\eta(x)$.
\end{theorem}

Having noted that transversals, complete mappings and orthomorphisms
are essentially the same thing, we will adopt the practice of expressing
our results in terms of transversals even when the original authors
used one of the other notions.

\medskip

As mentioned, this section is devoted to the case when our latin square
is $L_G$, the Cayley table of a finite group $G$.  The extra structure
in this case allows for much stronger results.  For example, suppose
we know of a transversal of $L_G$ that comprises a choice from each
row $i$ of an element $g_i$. Let $g$ be any fixed element of $G$. Then
if we select from each row $i$ the element $g_ig$ this will give a new
transversal and as $g$ ranges over $G$ the transversals so produced
will be mutually disjoint. Hence

\begin{theorem}\label{t:1all}
If $L_G$ has a single transversal then it has a decomposition into
disjoint transversals.
\end{theorem}

We saw in \sref{s:intro} that the question of which latin squares have
transversals has not been settled. The same is true for group tables, but
we are getting much closer to answering the question, building on
the pioneering work of Hall and Paige.

Consider the following five propositions:
\begin{enumerate}
\item[(i)]	$L_G$ has a transversal.
\item[(ii)] 	$L_G$ can be decomposed into disjoint transversals.
\item[(iii)] 	There exists a latin square orthogonal to~$L_G$.
\item[(iv)] 	There is some ordering of the elements of $G$, say $a_1,a_2,\dots,a_n$, such that {$a_1a_2\cdots a_n=\id$}, where $\id$ denotes the identity element of~$G$.
\item[(v)]	The Sylow 2-subgroups of $G$ are trivial or non-cyclic.
\end{enumerate}

The fact that (i), (ii) and (iii) are equivalent comes directly from
\tref{t:transMOLS} and \tref{t:1all}. Paige \cite{Paige}
showed that (i) implies (iv). Hall and Paige \cite{HP} then showed
that (iv) implies (v). They also showed that (v) implies (i) if $G$
is a soluble, symmetric or alternating group. They conjectured that
(v) is equivalent to (i) for all groups.

It was subsequently noted in \cite{DKIII} that both (iv) and (v) hold
for all non-soluble groups, which proved that (iv) and (v) are
equivalent.  A much more direct and elementary proof of this fact was
given in \cite{vlw03}.

To summarise:
\begin{theorem} \label{t:hallpaige}
(i)$\Leftrightarrow$(ii)$\Leftrightarrow$(iii)$\Rightarrow$(iv)$\Leftrightarrow$(v)
\end{theorem}

\begin{conjecture} \label{c:hallpaige}
(i)$\Leftrightarrow$(ii)$\Leftrightarrow$(iii)$\Leftrightarrow$(iv)$\Leftrightarrow$(v)
\end{conjecture}

As mentioned above, \cref{c:hallpaige} is known to be true for all
soluble, symmetric and alternating groups. It has also been shown for
many other groups including the linear groups $GL(2,q)$, $SL(2,q)$,
$PGL(2,q)$ and $PSL(2,q)$ (see \cite{EvansSL} and the references
therein).

After decades of incremental progress on \cref{c:hallpaige} there has
recently been what would appear to be a very significant
breakthrough. In a preprint Wilcox \cite{Wilcox} has claimed to reduce
the problem to showing it for the sporadic simple groups (of which the
Mathieu groups have already been handled in \cite{DVG}).  See
\cite{DKI}, \cite{EvansHP} or \cite{Wilcox} for further reading and
references on the Hall-Paige conjecture.

\bigskip

An immediate corollary of the proof of \tref{t:1all} is that for any
$G$ the number of transversals through a given entry of $L_G$ is
independent of the entry chosen.  Hence (see Theorem~3.5
of \cite{DKII}) we get:

\begin{theorem} 
The number of transversals in $L_G$ is divisible by $|G|$, the order of $G$.
\end{theorem}

McKay et al.\ \cite{MMW} also showed the following
simple results, in the spirit of Theorem~\ref{t:balasub}:

\begin{theorem} \label{t:oddsym}
The number of transversals in any symmetric latin square of order $n$
is congruent to $n$ modulo $2$.
\end{theorem}

\begin{corollary} \label{C:grpmod2}
Let $G$ be a group of order~$n$. If $G$ is abelian or $n$ is even then
the number of transversals in $G$ is congruent to $n$~modulo~$2$.
\end{corollary}

Corollary \ref{C:grpmod2} cannot be generalised to non-abelian groups
of odd order, given that the non-abelian group of order~21 has
826814671200 transversals.

\begin{theorem} \label{t:mod3}
If $G$ is a group of order $n\not\equiv1\mod3$ then the number of
transversals in $G$ is divisible by~$3$.
\end{theorem}

We will see below that the cyclic groups of small orders
$n\equiv1\mod3$ have a number of transversals which is not a multiple
of three. 

The semiqueens problem in \sref{s:enum} led to an investigation of
$z_n$, the number of transversals in the cyclic group of order $n$.
Let $z'_n = z_n/n$ denote the number of transversals through any given
entry of the cyclic square of order~$n$.  Since $z_n=z'_n=0$ for all
even $n$ by Theorem~\ref{t:hallpaige} we shall assume for the
following discussion that $n$ is odd.

The initial values of $z'_n$ are known from \cite{ys05} and
\cite{ys06}. They are $z'_1=z'_3=1$, $z'_5=3$, $z'_7=19$, $z'_9=225$,
$z'_{11}=3441$, $z'_{13}=79259$, $z'_{15}=2424195$, $z'_{17}=94471089$,
$z'_{19}=4613520889$, $z'_{21}=275148653115$, $z'_{23}=19686730313955$
and $z'_{25}=1664382756757625$. Interestingly, if we take these numbers
modulo~8 we find that this sequence begins
{1,1,3,3,1,1,3,3,1,1,3,3,1}. We know from Theorem~\ref{t:oddsym} that
$z'_n$ is always odd for odd $n$, but it is an open question
whether there is any deeper pattern modulo~4 or~8. We also know from
Theorem~\ref{t:mod3} that $z'_n$ is divisible by 3 when
$n\equiv2\mod3$. The initial terms of $\{z'_n\mod 3\}$ are
{1,1,0,1,0,0,2,0,0,1,0,0,2}.

An interesting fact about $z_n$ is that it is the number of {\em
diagonally cyclic latin squares} of order $n$ (in other words, the
number of quasigroups on the set $\{1,2,\dots,n\}$ which have the
transitive automorphism $(123\cdots n)$).  See \cite{DCLS} for a survey
on such objects.

\begin{table}[htb]
\centering
\begin{tabular}{|r|l|} \hline
	$n$& Number of transversals in groups of order $n$ \\
	\hline
	3	& 3\\
	4	& 0, 8\\
	5	& 15\\
	7	& 133\\
	8	& 0, 384, 384, 384, 384\\
	9	& 2025, 2241\\
	11& 37851\\
	12& 0, 198144, 76032, 46080, 0\\
	13& 1030367\\
	15& 36362925\\
	16& 0, 235765760, 237010944, 238190592, 244744192, 125599744,\\ 
		& 121143296, 123371520, 123895808, 122191872, 121733120,\\
		& 62881792, 62619648, 62357504\\
	17& 1606008513\\
	19& 87656896891\\ 
	20& 0, 697292390400, 140866560000, 0, 0\\ 
	21& 5778121715415, 826814671200\\ 
	23& 452794797220965\\ \hline
\end{tabular}
\caption{\label{Tab:groups}Transversals in groups of order $n\leq23$.}
\end{table}

We now discuss the number of transversals in general groups of small
order. For groups of order $n\equiv2\mod4$ there can be no
transversals, by \tref{t:hallpaige}. For each other order~$n\leq23$ the
number of transversals in each group is given in
Table~\ref{Tab:groups}. The groups are ordered according to the
catalogue of Thomas and Wood~\cite{tw80}.  The numbers of transversals
in abelian groups of order at most~$16$ and cyclic groups of order at
most $21$ were obtained by Shieh et al~\cite{shh00}. The remaining
values in Table \ref{Tab:groups} were computed by
Shieh~\cite{ys05}. McKay et al.\ \cite{MMW} then independently
confirmed all counts except those for cyclic groups of order $\ge21$,
correcting one misprint in Shieh~\cite{ys05}.

Bedford and Whitaker~\cite{bw99} offer an explanation for why all the
non-cyclic groups of order~8 have 384 transversals. The groups of
order 4, 9 and 16 with the most transversals are the elementary
abelian groups of those orders. Similarly, for orders 12, 20 and~21
the group with the most transversals is the direct sum of cyclic
groups of prime order.  It is an open question whether such a
statement generalises to all $n$.

By Corollary~\ref{C:grpmod2} we know that in each case covered by
Table~\ref{Tab:groups} (except the non-abelian group of order 21), the
number of transversals must have the same parity as the order of the
square. It is remarkable though, that the groups of even order have a
number of transversals which is divisible by a high power of
2. Indeed, any 2-group of order $n\leq16$ has a number of transversals
which is divisible by~$2^{n-1}$. It would be interesting to know
if this is true for general~$n$.

\section*{\centerline{5. Generalised transversals}}%\label{s:plex}

There are several ways to generalise the notion of a transversal.
We have already seen one of them, namely the partial transversals in
\sref{s:partial}. In this section we collect results on another
generalisation, namely plexes.

A \emph{$k$-plex} in a latin square of order $n$ is a set of $kn$
entries which includes $k$ representatives from each row and each
column and of each symbol. A transversal is a $1$-plex. 
The marked entries form
a $3$-plex in the following square:

\begin{equation}\label{e:3plex}
\begin{matrix}
1^*&2\phantom{^*}&3\phantom{^*}&4^*&5\phantom{^*}&6^*\cr
2^*&1\phantom{^*}&4\phantom{^*}&3^*&6^*&5\phantom{^*}\cr
3\phantom{^*}&5^*&1\phantom{^*}&6\phantom{^*}&2^*&4^*\cr
4\phantom{^*}&6\phantom{^*}&2^*&5\phantom{^*}&3^*&1^*\cr
5^*&4^*&6^*&2\phantom{^*}&1\phantom{^*}&3\phantom{^*}\cr
6\phantom{^*}&3^*&5^*&1^*&4\phantom{^*}&2\phantom{^*}\cr
\end{matrix}
\end{equation}
The name $k$-plex was coined in \cite{iw02} only recently. It is a natural
extension of the names duplex, triplex, and quadruplex which have been
in use for many years (principally in the statistical literature, such
as \cite{BFinnI}) for 2, 3 and 4-plexes.

The entries not included in a $k$-plex of a latin square $L$ of order
$n$ form an $(n-k)$-plex of $L$. Together the $k$-plex and its
complementary $(n-k)$-plex are an example of what is called an
\emph{orthogonal partition} of $L$. For discussion of orthogonal partitions
in a general setting see Gilliland \cite{BGill} and
Bailey~\cite{BBail}.  For our purposes, if $L$ is decomposed into
disjoint parts $K_1$, $K_2,\dots,K_d$ where $K_i$ is a $k_i$-plex then
we call this a $(k_1,k_2,\dots,k_d)$-partition of $L$. A case of
particular interest is when all parts are the same size, $k$. We call
such a partition a {\it $k$-partition\/}. For example, the marked
$3$-plex and its complement form a $3$-partition of the square in
\eref{e:3plex}.  By \tref{t:transMOLS}, finding a 1-partition of a square is
equivalent to finding an orthogonal mate.

Some results about transversals generalise directly to other plexes,
while others seem to have no analogue. 
\tref{t:balasub} and \tref{t:1all} 
seem to be in the latter class, as observed in
\cite{MMW} and \cite{iw02} respectively.  However, Theorems
\ref{Tqstep} and \ref{t:hallpaige} showed that not every square has a
transversal, and exactly the same arguments work for any $k$-plex
where $k$ is odd \cite{iw02}.

\begin{theorem}\label{Tqstepplex}
Suppose that $q$ and $k$ are odd integers and $m$ is even.
No $q$-step type latin square of order $mq$ possesses a $k$-plex.
\end{theorem}

\begin{theorem}\label{Tgpqstepplex}
Let $G$ be a group of finite order $n$ with a non-trivial 
cyclic Sylow $2$-subgroup. The Cayley table of $G$ contains no $k$-plex
for any odd $k$ but has a $2$-partition and hence contains a $k$-plex for
every even $k$ in the range $0\leq k\leq n$.
\end{theorem}

The situation for even $k$ is quite different to the odd case.  Rodney
\cite[p.105]{BColD} conjectures that every latin square has a
duplex. This conjecture was strengthened in \cite{iw02} to the
following:

\begin{conjecture}\label{c:maxdupe}
Every latin square has the maximum possible number of disjoint
duplexes.  In particular, every latin square of even order has a
$2$-partition and every latin square of odd order has a
$(2,2,2,\dots,2,1)$-partition.
\end{conjecture}

Note that this conjecture also strengthens \cref{c:ryser}. It also
implies that every latin square has $k$-plexes for every even
value of $k$ up to the order of the square.

\cref{c:maxdupe} is true for all latin squares of orders $\le8$ and
for all soluble groups (see \cite{vlw03,iw02}).
Depending on whether a soluble group has a non-trivial cyclic
Sylow $2$-subgroup, it either has a $k$-plex for all possible $k$, or has
them for all possible even $k$ but no odd $k$. If the Hall-Paige
conjecture could be proved it would completely resolve the existence
question of plexes in groups, and these would remain the only two
possibilities.  It is worth noting that other scenarios occur for
latin squares which are not based on groups. For example, the square
in \eref{e:3plex} has no transversal but clearly does have a $3$-plex.
It is conjectured in \cite{iw02} that there exist arbitrarily large
latin squares of this type.

\begin{conjecture}\label{cjthrnotran}
For all even $n>4$ there exists a latin square of 
order $n$ which has no transversal but does contain a $3$-plex. 
\end{conjecture}

Another possibility was shown by a family of squares constructed
in \cite{egan1}.

\begin{theorem}
For all even $n>2$ there exists a latin square of order $n$ which has
$k$-plexes for every odd value of $k$ between $\lfloor n/4\rfloor$ and
$\lceil 3n/4\rceil$, but not for any odd value of $k$ outside this
range.
\end{theorem}

Interestingly, there is no known example of odd integers $a<b<c$ and a
latin square which has an $a$-plex and a $c$-plex but no $b$-plex.

\medskip

%Following Finney~\cite{BFinnIV}, we say that two plexes in the same
%square are {\it parallel\/} if they have no filled cells in common.

The union of an $a$-plex and a disjoint $b$-plex of a latin square $L$
is an $(a+b)$-plex of $L$. However, it is not always possible to
split an $(a+b)$-plex into an $a$-plex and a disjoint $b$-plex.
Consider a duplex which consists of $\frac12n$ disjoint intercalates
(latin subsquares of order 2).  Such a duplex does not contain a
partial transversal of length more than $\frac12n$, so it is a long
way from containing a 1-plex.

We say that a $k$-plex is {\it indivisible\/} if it contains no
$c$-plex for $0<c<k$. The duplex just described is
indivisible. Indeed, for every $k$ there is a indivisible $k$-plex in
some sufficiently large latin square. This was first shown in
\cite{iw02}, but ``sufficiently large'' in that case meant quadratic
in $k$. This was improved to linear in \cite{egan2} as a corollary of
the following result.

\begin{theorem}\label{t:indivhalf}
For every $k\ge2$ there exists a latin square of order $2k$ which contains
two disjoint indivisible $k$-plexes.
\end{theorem}

\tref{t:indivhalf} means that some squares can be split in ``half'' in
a way that makes no further division possible. Experience with latin
squares suggests that they generally have a vast multitude of
partitions into various plexes, which in one sense means that latin
squares tend to be a long way from being indivisible.  This makes
\tref{t:indivhalf} slightly surprising. 

It is a wide open question for what values of $k$ and $n$ there is a
latin square of order $n$ containing an indivisible $k$-plex.
However, Bryant et al.~\cite{egan2} found the answer when $k$ is small
relative to $n$.

\begin{theorem}\label{t:indiv5th}
Let $n$ and $k$ be positive integers satisfying $5k\le n$. Then there exists
a latin square of order $n$ containing an indivisible $k$-plex.
\end{theorem}

\bigskip

So far we have essentially looked at questions where we start with a
latin square and ask what sort of plexes it might have.  To complete
the section we consider the reverse question. We want to start with a
plex and ask what latin squares it might be contained in. Strictly
speaking this is a silly question, since we defined a plex in terms of
its host latin square, which therefore is the only possible answer.
However, suppose we define a \emph{$k$-homogeneous partial latin
square} of order $n$ to be an $n\times n$ array in which each cell is
either blank or filled (the latter meaning that it contains one of
$\{1,2,\dots,n\}$), and which has the properties that (i)~no symbol
occurs twice within any row or column, (ii)~each symbol occurs $k$
times in the array, (iii)~each row and column contains exactly $k$
filled cells. (The standard definition of a homogeneous partial latin
square is slightly more general. However, once empty rows and columns
have been deleted, it agrees with ours.) We can then sensibly ask
whether this $k$-homogeneous partial latin square is a $k$-plex.  If
it is then we say the partial latin square is \emph{completable}
because the blank entries can be filled in to produce a latin square.

\begin{theorem}\label{Tnoncomp}
If $1<k<n$ and $k>\frac14n$ then there exists a $k$-homogeneous partial
latin square of order $n$ which is not completable.
\end{theorem}

Burton \cite{Burton}, and Daykin and H\"aggkvist \cite{BDyHg}
independently conjecture that if $k\leq\frac14n$ then every $k$-plex
is completable. It seems certain that for $k$ sufficiently small
relative to $n$, every $k$-plex is completable. This has already been
proved when $n\equiv0\mod16$ in \cite{BDyHg}.  The following partial
extension result due to Burton \cite{Burton} also seems relevant.

\begin{theorem}\label{thpartext}
For $k\leq\frac14n$ every $k$-homogeneous partial latin square of order $n$
is contained in a $(k+1)$-homogeneous partial latin square of order~$n$.
\end{theorem}

%In the statistical literature (eg. Finney \cite{BFinnI}
%\cite{BFinnIV}) a transversal is sometimes called a directrix;

\section*{\centerline{6. Covering radii for sets of permutations}}%\label{s:kezsnev}

\newcommand{\crad}{\mathop{\mathrm{cr}}}

A novel approach to \cref{c:ryser} and \cref{c:brualdi}
has recently been opened up by Andre K\'ezdy and Hunter Snevily.
To explain this interesting new approach, we need to introduce 
some terminology.

Consider the \emph{symmetric group} $S_n$ as a metric space equipped
with \emph{Hamming distance}. That is, the distance between two
permutations $g,h\in S_n$ is the number of points at which they
disagree ($n$ minus the number of fixed points of $gh^{-1}$). Let $P$
be a subset of $S_n$.  The \emph{covering radius} $\crad(P)$ of $P$ is
the smallest $r$ such that the balls of radius $r$ with centres at the
elements of $P$ cover the whole of $S_n$.  In other words every
permutation is within distance $r$ of some member of $P$, and $r$ is
chosen to be minimal with this property.

\begin{theorem}
Let $P\subseteq S_n$ be a set of permutations. If $|P|\le n/2$, then
$\crad(P)=n$.  However, there exists $P$ with $|P|=\lfloor n/2\rfloor+1$
and $\crad(P)<n$.
\label{p1}
\end{theorem}

This result raises an obvious question. Given $n$ and $s$, what is the
smallest $m$ such that there is a set $S$ of permutations with $|S|=m$
and $\crad(S)\le n-s$? We let $f(n,s)$ denote this minimum value $m$.
This problem can also be interpreted in graph-theoretic
language. Define the graph $G_{n,s}$ on the vertex set $S_n$, with two
permutations being adjacent if they agree in at least $s$ places.  Now
the size of the smallest dominating set in $G_{n,s}$ is $f(n,s)$.

Theorem~\ref{p1} shows that
$f(n,1)=\lfloor n/2\rfloor+1$. Since any two distinct permutations have
distance at least~$2$, we see that $f(n,n-1)=n!$ for $n\ge2$. Moreover,
$f(n,s)$ is a monotonic increasing function of $s$ (by definition).

The next case to consider is $f(n,2)$.  K\'ezdy and Snevily made the
following conjecture in unpublished notes.

\begin{conjecture}\label{c:KS}
If $n$ is even, then $f(n,2)=n$; if $n$ is odd, then $f(n,2)>n$.
\end{conjecture}

The K\'ezdy--Snevily conjecture has several connections with 
transversals. The rows of a latin square of order~$n$ form a \emph{sharply
transitive set} of permutations (that is, exactly one permutation
carries $i$ to $j$, for any $i$ and $j$); and every sharply transitive
set is the set of rows of a latin square.

\begin{theorem}
Let $S$ be a sharply transitive subset of~$S_n$. Then $S$ has covering
radius at most $n-1$, with equality if and only if the corresponding latin
square has a transversal.
\end{theorem}

\begin{corollary}\label{c:transv}
If there exists a latin square of order~$n$ with no transversal, then
$f(n,2)\le n$. In particular, this holds for $n$ even.
\label{onetwo}
\end{corollary}

Hence \cref{c:KS} implies \cref{c:ryser},
as K\'ezdy and Snevily observed.  In fact a stronger result holds:

\begin{theorem}\label{p:n-2}
If $S$ is the set of rows of a latin square $L$ of order~$n$
with no transversal, then $S$ has covering radius $n-2$.
\end{theorem}

The following result is due to K\'ezdy and Snevily. See \cite{CW} for a
proof.

\begin{theorem}
\cref{c:KS} implies \cref{c:brualdi}.
\end{theorem}

In other words, to solve the longstanding Ryser and Brualdi conjectures
it may suffice to answer this: How small can we make a subset $S\subset S_n$
which has the property that every permutation in $S_n$ agrees with some
member of $S$ in at least two places?

In Corollary~\ref{c:transv} we used latin squares to find an upper
bound for $f(n,2)$ when $n$ is even.  For odd $n$ we can also find
upper bounds based on latin squares.  The idea is to choose a latin
square with few transversals, or whose transversals have a particular
structure, and add a small set of permutations meeting each
transversal twice. For $n=5,7,9$, we now give a latin square for
which a single  extra permutation suffices, showing that $f(n,2)\le n+1$ in
these cases.
{\small
\[
\begin{matrix}
1&2&3&4&5\\
2&1&4&5&3\\
3&5&1&2&4\\
4&3&5&1&2\\
5&4&2&3&1\\[0.5ex]
1&3&4&2&5
\end{matrix}\qquad\quad
\begin{matrix}
1&2&3&4&5&6&7\\
2&3&1&5&4&7&6\\
3&1&2&6&7&4&5\\
4&5&6&7&1&2&3\\
5&4&7&1&6&3&2\\
6&7&4&2&3&5&1\\
7&6&5&3&2&1&4\\[0.5ex]
3&2&1&7&6&5&4
\end{matrix}\qquad\quad
\begin{matrix}
1&3&2&4&6&5&7&9&8\\
2&1&3&5&4&6&8&7&9\\
3&2&1&7&9&8&4&6&5\\
4&6&5&9&8&7&1&3&2\\
5&4&6&8&7&9&3&2&1\\
6&5&4&2&1&3&9&8&7\\
7&9&8&1&3&2&5&4&6\\
8&7&9&3&2&1&6&5&4\\
9&8&7&6&5&4&2&1&3\\[0.5ex]
5&4&6&1&3&2&9&8&7
\end{matrix}
\]
}

In general, we have the following:

\begin{theorem} $f(n,2)\le\frac43n+O(1)$ for all $n$.
\end{theorem}

The reader is encouraged to seek out \cite{CW}
and the survey by Quistorff \cite{Quistorff} for more
information on covering radii for sets of permutations.

\section*{\centerline{7. Concluding Remarks}}%\label{s:conclusion}

We have only been able to give the briefest of overviews of the
fascinating subject of transversals in this survey. Space constraints
have forced the omission of much worthy material, including proofs
of the theorems quoted.  However, even this brief skim across the
surface has shown that many basic questions remain unanswered and much
work remains to be done. 

%The subject is littered with tantalising
%conjectures.  Even the theorems in many cases seem to be far from best
%possible, leaving openings for future improvements. It is hoped that
%this survey may in a small way aid in achieving such improvements.

  \let\oldthebibliography=\thebibliography
  \let\endoldthebibliography=\endthebibliography
  \renewenvironment{thebibliography}[1]{%
    \begin{oldthebibliography}{#1}%
      \setlength{\parskip}{0.3ex plus 0.1ex minus 0.1ex}%
      \setlength{\itemsep}{0.3ex plus 0.1ex minus 0.1ex}%
  }%
  {%
    \end{oldthebibliography}%
  }

\noindent
\footnotesize{School of Mathematical Sciences\\[-0.5ex]
Monash University\\[-0.5ex]
Vic~3800, Australia\\
e-mail: ian.wanless@sci.monash.edu.au }

\end{document}